\documentclass[12pt]{article}
\RequirePackage[colorlinks,citecolor=blue,urlcolor=blue,linkcolor=blue]{hyperref}
\hypersetup{
colorlinks = true,
citecolor=blue,
urlcolor=blue,
linkcolor=blue,
pdfauthor = {Alexey Kuznetsov},
pdfkeywords = {L-function, completely multiplicative function, Euler product, functional equation, subordinator, Kendall's identity},
pdftitle = {On Dirichlet series and functional equations},
pdfpagemode = UseNone
}
\usepackage{graphicx,xspace,colortbl}
\usepackage{amsmath,amsthm,amsfonts,amssymb}
\usepackage{color}
\usepackage{enumerate}
\usepackage{fancybox}
\usepackage{epsfig}
\usepackage{subfig}
\usepackage{pdfsync}
    \def\qed{\hfill$\sqcap\kern-8.0pt\hbox{$\sqcup$}$\\}
    \def\beq{\begin{eqnarray}}
    \def\eeq{\end{eqnarray}}
    \def\beqq{\begin{eqnarray*}}
    \def\eeqq{\end{eqnarray*}}

\DeclareMathOperator{\re}{Re}
    \def\p{{\mathbb P}}
    \def\e{{\mathbb E}}
    \def\r{{\mathbb R}}
    \def\C{{\mathbb C}} 
    	
    \def\d{{\textnormal d}}
    \def\i{{\textnormal i}}

    \def\ee{{\textnormal e}}

\newtheorem{theorem}{Theorem}

\newtheorem{corollary}{Corollary}
\theoremstyle{definition}

\newtheorem{remark}{Remark}

\title{On Dirichlet series and functional equations}
\author{
{Alexey Kuznetsov
\footnote{Dept. of Mathematics and Statistics,  York University,
4700 Keele Street, Toronto, ON, M3J 1P3, Canada.  \newline
E-mail:  kuznetsov@mathstat.yorku.ca  } 
 }}
 
 \date{\today}

\begin{document}
\maketitle

\begin{abstract} 
There exist many explicit evaluations of Dirichlet series. Most of them are constructed via the same approach: by taking products or powers of  Dirichlet series with a known Euler product representation. In this paper we derive a result of a new flavour: we give the Dirichlet series representation to solution $f=f(s,w)$ of the functional equation 
$L(s-wf)=\exp(f)$, where $L(s)$ is the L-function corresponding to a completely multiplicative function. Our result seems to be a Dirichlet series analogue of the well known Lagrange-B\"urmann formula for power series. The proof is probabilistic in nature and is based on Kendall's identity, which arises in the fluctuation theory of L\'evy processes. 
\end{abstract}

{\vskip 0.15cm}

\section{Introduction and the main result}

Let $a  : {\mathbb N} \mapsto \C$ be a completely multiplicative function, that is, $a(mn)=a(m)a(n)$ for all 
$m, n \in {\mathbb N}$. Denote by $L(s)$ the  corresponding 
L-function
\begin{equation}\label{Dirichlet_series_for_L}
L(s)=\sum\limits_{n\ge 1} \frac{a(n)}{n^s}.
\end{equation}
We assume that the above series converges absolutely for $\re(s)\ge \sigma$.
Complete multiplicativity of $a(n)$ implies that $L(s)$ can be expressed as an absolutely convergent Euler product
$$
L(s)=\prod\limits_{p} (1-a(p)p^{-s})^{-1}, \;\;\; \re(s)\ge \sigma, 
$$
where the product is taken over all prime numbers $p$. 
The following functions will play the key role in what follows: for $n \in {\mathbb N}$ and $z\in \C$ we define
\begin{equation}\label{def_tau_nt}
d_z(n):=\prod\limits_{p^j | n} \binom{j+z-1}{j},  \qquad
\tilde d_{z}(n):=z^{-1} d_z(n). 
\end{equation}
The function $d_z(n)$ is multiplicative and it is called {\it the general divisor function}, see 
\cite{Ivic}[Section 14.6]. 
Starting from the Euler product representation for $\zeta(s)$ and writing the terms 
$(1-p^{-s})^{-z}$ as binomial series in $p^{-s}$, it is easy to see that 
\begin{equation}\label{eqn_power_zeta}
\zeta(s)^z=\sum\limits_{n\ge 1} \frac{ d_z(n)}{n^s}, \;\;\; z\in \C, \;\;\; s>1.
\end{equation}

The multiplicative function $d_z(n)$ is well known in the literature. Selberg \cite{Selberg} has obtained the main term of the asymptotics of $D_z(x):=\sum_{n\le x} d_z(n)$ as $x\to +\infty$; the information about higher-order terms can be found in \cite{Ivic}[Theorem 14.9].    The function $d_k(n)$ (for integer $k\ge 2$) is known simply as {\it the divisor function}
(see \cite{Hardy_Wright}[Section 17.8] or \cite{gould2008}). The name comes from the following fact
$$
d_k(n)=\sum\limits_{\substack{m_1 m_2 \cdots m_k=n \\ m_i\ge 1}} 1,
$$
which follows from from \eqref{eqn_power_zeta}. 
In other words, $d_k(n)$ counts the number of ways of expressing $n$ as an ordered product of $k$ positive factors (of which any number may be unity). For example, $d_2(n)$ is the number of divisors of $n$, which is commonly denoted by $d(n)$. Also, note that  for all $n\ge 2$ the function 
$z\mapsto  \tilde d_{z}(n)$ is a polynomial of degree $\Omega(n)-1$, where $\Omega(n)$ is the total number of prime factors of $n$. In particular, $\tilde d_{z}(n) \equiv 1$ if and only if $n$ is a prime number.

Let us denote 
\begin{equation}\label{def_D_sigma_R}
D_{\sigma,\rho}:=\{(s,w) \in \C^2 \; : \; \re(s)\ge \sigma, \; |w| \le \rho\}.
\end{equation}
The following theorem is our main result.

\begin{theorem}\label{thm1} Assume that   the Dirichlet series \eqref{Dirichlet_series_for_L}, which corresponds to 
a completely multiplicative function $\{a(n)\}_{n \in {\mathbb N}}$,  converges absolutely 
for $\re(s) \ge \sigma$. Denote 
\begin{equation}\label{def_gamma}
\gamma:=\ln \Big(\sum\limits_{n\ge 1} \frac{|a(n)|}{n^{\sigma}}\Big). 
\end{equation}
Then for any $\rho>0$:
\begin{itemize}
\item[(i)] The series
\begin{equation}\label{def_f_s_w}
f(s,w):=\sum\limits_{n\ge 2}  \tilde d_{w\ln(n)}(n) \times \frac{a(n)}{n^s} 
\end{equation}
converges absolutely and uniformly in $(s,w) \in D_{\sigma+\gamma\rho,\rho}$ and satisfies $|f(s,w)|<\gamma$ in this region; 
\item[(ii)] 
The function $f(s,w)$ solves the functional equation
\begin{equation}\label{inverse_eqn}
L(s-wf(s,w))=\exp(f(s,w)), \;\;\; (s,w)  \in D_{\sigma+\gamma\rho,\rho}; 
\end{equation}
\item[(iii)]
For any $v \in \C$ and $(s,w) \in D_{\sigma+\gamma\rho,\rho}$ the following identity is true
\begin{equation}\label{eqn_exp_v_f_s_w_}
1+v\sum\limits_{n\ge 2}\tilde d_{v+w\ln(n)}(n) \times \frac{a(n)}{n^s}= \exp(vf(s,w)).
\end{equation}
\end{itemize}
\end{theorem}

The proof of Theorem \ref{thm1} is presented in the next section. 

\begin{remark}
Let us consider what happens with formulas \eqref{def_f_s_w} and \eqref{eqn_exp_v_f_s_w_} when  $w=0$. Note that $\tilde d_0(n)=1/j$ if $n=p^j$ for some prime $p$ and $\tilde d_0(n)=0$ otherwise. Then we can write $\tilde d_0(n) = \Lambda(n)/\ln(n)$, where 
 the von Mangoldt function $\{\Lambda(n)\}_{n\in {\mathbb N}}$ is defined as follows
\begin{align}\label{def_Lambda}
\Lambda(n)=
\begin{cases} 
\ln(p), \;\;\; {\textnormal{ if $n=p^k$ for some prime $p$ and integer $k\ge 1$}}, \\
0, \;\;\qquad  {\textnormal{ otherwise}}.
\end{cases}
\end{align}
Using the above result and
\eqref{def_f_s_w} we obtain
\begin{equation}\label{eqn_f_s_0}
f(s,0)=\sum\limits_{n\ge 2} \frac{\Lambda(n)a(n) }{\log(n)n^s }=\ln(L(s)) 
\end{equation}
for $\re(s)\ge \sigma$. 
Formula \eqref{eqn_f_s_0} confirms the functional identity 
\eqref{inverse_eqn} in the case $w=0$. Formula \eqref{eqn_exp_v_f_s_w_} in the case $w=0$ also becomes a trivial identity
$L(s)^v=\exp(v \ln L(s))$ (see equation \eqref{eqn_L_s_t} below). 
\end{remark}

\begin{remark}
Formula \eqref{def_f_s_w}, which gives a solution to the functional equation \eqref{inverse_eqn}, has some similarities to the Lagrange-B\"urmann inversion formula for analytic functions. Let us remind what the latter result states. Consider a function $\psi$, which is analytic in a neighbourhood of $w=0$ and satisfies $\psi(0)\neq 0$. Let $w=g(z)$ denote the solution of  $z\psi(w)=w$. Then $g(z)$ can be represented as a convergent Taylor series 
\begin{equation}\label{eqn_Burmann}
g(z)=\sum\limits_{n\ge 1} \lim\limits_{w\to 0} \Big[ \frac{\d^{n-1}}{\d w^{n-1}} \psi(w)^n \Big] \frac{z^n}{n!},
\end{equation}
which converges in some neighbourhood of $z=0$. 
Note that both formulas \eqref{eqn_Burmann} and \eqref{def_f_s_w} are based on the coefficients of the expansion of a power of the original function in a certain basis (the basis consists of power functions $z^n$ in the case of the Lagrange-B\"urmann inversion formula and exponential functions $n^{-s}$ in the case of formula \eqref{def_f_s_w}). 
\end{remark}

\begin{remark}
Formula \eqref{eqn_power_zeta} implies the well-known result
\begin{equation}\label{semigroup_result1}
d_{t+s}(n)=\sum\limits_{k|n} d_{t}(k) \times d_s(n/k), \;\;\; t,s \in \C.  
\end{equation}
Similarly, formula \eqref{eqn_exp_v_f_s_w_} implies the following more general result
\begin{equation}\label{semigroup_result}
(t+s) \tilde d_{t+s+w\ln(n)}(n) =
ts\sum\limits_{ k|n}
\tilde d_{t+w\ln(k)}(k)  \times \tilde d_{s+w\ln(n/k)}(n/k), 
\;\;\; t,s,w \in \C. 
\end{equation}
Note that \eqref{semigroup_result1} is a special case of  \eqref{semigroup_result} with $w=0$ and that 
both sides of \eqref{semigroup_result} are polynomials in variables $(t,s,w)$. It would be an interesting exercise to try to find an elementary proof of \eqref{semigroup_result}. 
\end{remark}

Next we present a corollary of Theorem \ref{thm1}, its proof is postponed until section \ref{section_cor_proofs}.  
Whenever we use $\ln(L(s))$ in what follows, we will always assume that $\re(s)\ge \sigma$ and the branch of logarithm is chosen so that 
\eqref{eqn_f_s_0} holds (another way to fix the branch of logarithm is to require that $\ln(L(s))\to 0$ as $s\to +\infty$). 

\begin{corollary}\label{cor1} 
Let $\gamma$ and $f(s,w)$ be defined as in \eqref{def_gamma} and \eqref{def_f_s_w}. 
Then for all $(s,w) \in \C^2$ satisfying 
$\re(s) \ge \sigma + 2 \gamma |w|$ we have
$f(s + w \ln(L(s)),w)=\ln(L(s))$. 
\end{corollary}

The above result can be used to obtain new explicit evaluations of infinite series. For example, consider the case when $a(n)=1$ for all $n$,
so that $L(s)=\zeta(s)$, which is the Riemann zeta-function.  Then, using
formula \eqref{eqn_exp_v_f_s_w_} and Corollary \ref{cor1} with $\sigma=1.4$, $s=2$ (so that $\ln(\zeta(2))=\ln(\pi^2/6)$) we obtain the following explicit result
\begin{equation}\label{explicit_series}
\sum\limits_{n\ge 2} \frac{\tilde d_{v+w\ln(n)}(n)}{n^z}=\frac{\big( \pi^2/6 \big)^v-1}{v},
\end{equation} 
which is valid for all $v\in \C$, $z$ in the disk $\{ z\in \C\; : \; |z-2|\le 0.13\}$, and $w=(z-2)/\ln(\pi^2/6)$. A curious fact is that the sum of the series \eqref{explicit_series} does not depend on $z$. As we have mentioned above, the functions 
$v\mapsto \tilde d_{v+ w \ln(n)}(n)$ are polynomials of degree $\Omega(n)-1$, thus formula \eqref{explicit_series} can be viewed as 
an expansion of the entire function $v \in \C \mapsto  \big(\big( \pi^2/6 \big)^v-1\big)/v$ in such an unusual polynomial basis.  There are two natural questions that arise from identity \eqref{explicit_series}: (i) What is the largest domain of $z$ for which the series converges absolutely or conditionally? (ii) Is it  possible to find an elementary proof of \eqref{explicit_series}?

\section{Proof of Theorem \ref{thm1}}\label{section_proofs}

The proof of Theorem \ref{thm1} will proceed in two stages. First we will prove Theorem \ref{thm1} in the special case when $v>0$, $w>0$ and $a(n)\ge 0$ for all $n \in {\mathbb N}$. 
This proof is probabilistic in nature and it is based on the theory of L\'evy processes, see \cite{Kyprianou}. In the second stage we will complete the proof of Theorem \ref{thm1} by generalizing our earlier result to complex values $v$, $w$ and $a(n)$ by an analytic continuation argument.

For convenience of the reader, we will first review several key facts from the theory of L\'evy processes, which will be required in our proof (one may wish to consult the books \cite{Bertoin} and \cite{Kyprianou} for more detailed information).
 A one-dimensional stochastic process $X=\{X_t\}_{t\ge 0}$ is called a {\it subordinator} if it has stationary and independent increments and if its paths (functions  $t\in [0,\infty) \mapsto X_t$) are increasing almost surely. We will always assume that $\p(X_0=0)=1$.
A probability measure $\nu(\d x)$ supported on $[0,\infty)$ is called {\it infinitely divisible} if for any $n=2,3,4,\dots$ there exists a probability measure
$\nu_n(\d x)$  such that $\nu=\nu_n * \nu_n * \cdots * \nu_n$ ($\nu$ is an $n$-fold convolution of the measure $\nu_n$). It is known that subordinators stand in one-to-one correspondence with infinitely divisible measures: for any subordinator $X$ the measure $\p(X_1\in \d x)$ is infinitely divisible and for any infinitely divisible measure $\nu$ supported on $[0,\infty)$ there exists a unique subordinator $X$ such that $\p(X_1\in \d x)=\nu(\d x)$. 

Let  $X$ be a subordinator and   $\ee(\kappa)$ be an exponential random variable with mean $1/\kappa$, 
independent of $X$. We define a new process via 
\begin{equation*}
\tilde X_t=\begin{cases}
X_t, \;\;\;\;\;\; {\textnormal{ if }} t<\ee(\kappa), \\
+\infty, \;\;\;\; {\textnormal{ if }} t\ge \ee(\kappa).
\end{cases}
\end{equation*}
The process $\tilde X$ is called a {\it killed subordinator}. Note that killed subordinators satisfy $\p(\tilde X_t \in [0,\infty))=\p(\ee(\kappa)>t)=\exp(-\kappa t)$, thus the measures $\p(\tilde X_t\in \d x)$ are sub-probability measures.

Any subordinator (including killed ones) can be described through an associated 
{\it Bernstein function}  $\phi_X(z)$ via the identity
\begin{equation}\label{Laplace_X_t}
\e\Big[e^{-z X_t}\Big]=\int_{[0,\infty)} e^{-z x } \p(X_t \in \d x)=e^{-t \phi_X(z)}, \;\;\; \re(z)\ge 0, \; t>0. 
\end{equation}
The above identity expresses the fact that the probability measures $\mu_t(\d x)=\p(X_t \in \d x)$ form a convolution semigroup on $[0,\infty)$, that is $\mu_t * \mu_s=\mu_{t+s}$. 
The L\'evy-Khintchine formula tells us that any Bernstein function has an integral representation 
\begin{equation}\label{Laplaceexponent}
\phi_X(z) = \kappa + \delta z + \int_{(0,\infty)} (1-{\rm e}^{-z x})\Pi(\d x), \;\;\; \re(z)\ge 0,
\end{equation}
for some  $\kappa\geq 0$, $\delta\geq 0$ and a positive measure $\Pi(\d x)$, supported on $(0,\infty)$, which satisfies the integrability condition $\int_{(0,\infty)}(1\wedge x)\Pi(\d x)<\infty$. The constant $\kappa$ is called {\it the killing rate}, $\delta$ is called {\it the linear drift coefficient} and the measure 
$\Pi(\d x)$ is called {\it the L\'evy measure}. The L\'evy measure describes the distribution of jumps of the process $X$. 
The killing rate and the drift can be recovered from the Bernstein function as follows:
$$
\kappa=\phi_X(0) {\textnormal{  and  }} \delta=\lim\limits_{z\to +\infty} \frac{\phi_X(z)}{z}. 
$$
See the excellent book \cite{Schilling} for more information on Bernstein functions. 

In this paper we will only be working with a rather simple class of subordinators -- the ones that have compound Poisson jumps. The L\'evy measure of such a process has finite  mass $\lambda=\Pi([0,\infty))<\infty$ and the process itself can be constructed as follows. Take a sequence of independent and identically distributed random variables $\xi_i$, having distribution $\p(\xi_i \in \d x)=\lambda^{-1} \Pi(\d x)$, and take an independent Poisson process $N=\{N_t\}_{t\ge 0}$ with intensity $\lambda$ (that is, $\e[N_t]=\lambda t$). Then the pathwise definition of the subordinator $X$ is 
$$
X_t=\delta t + \sum\limits_{i=1}^{N_t} \xi_i, \;\;\; t\ge 0. 
$$

\begin{figure}
\centering
\subfloat[][]{\label{f11}\includegraphics[height =5.5cm]{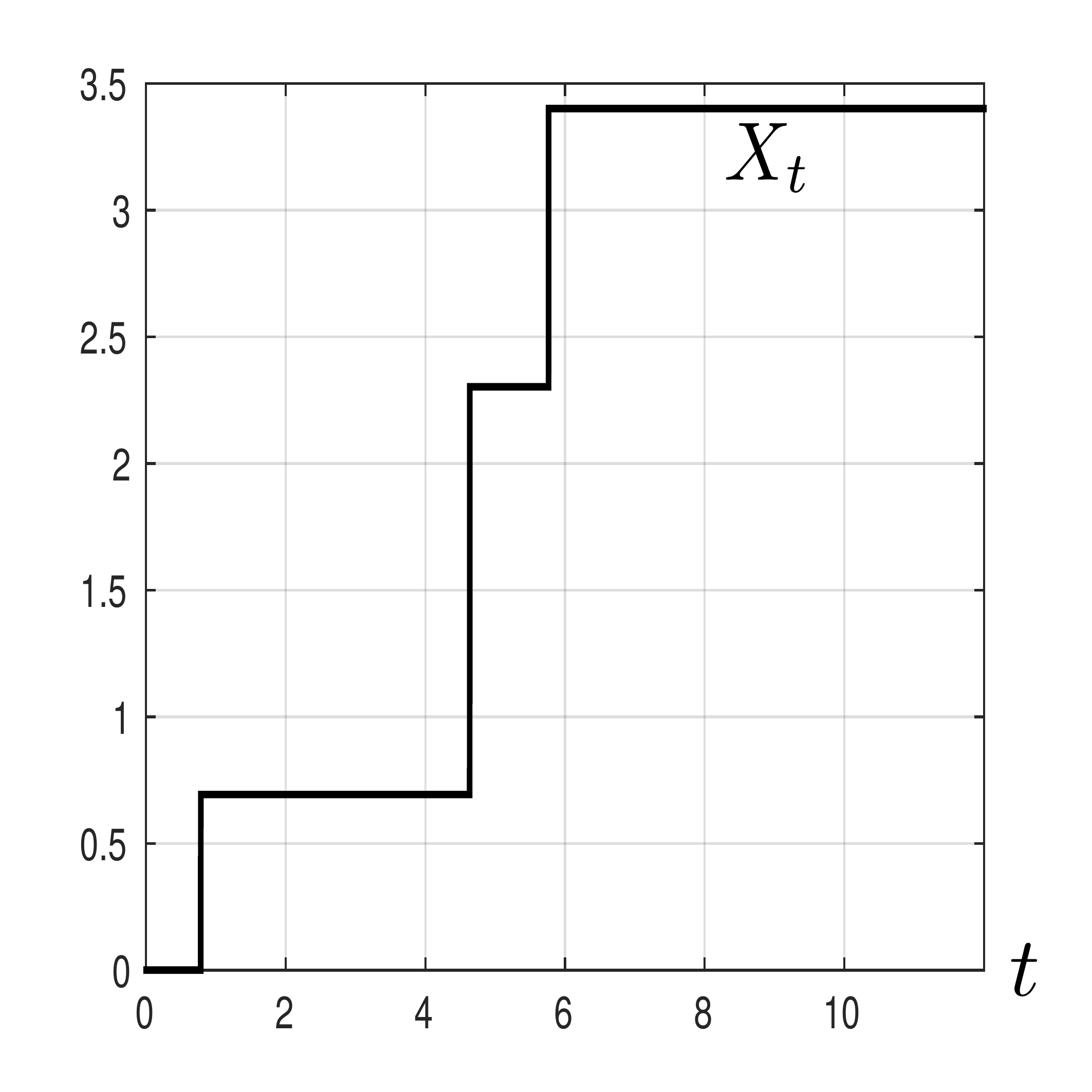}}
\subfloat[][]{\label{f12}\includegraphics[height =5.5cm]{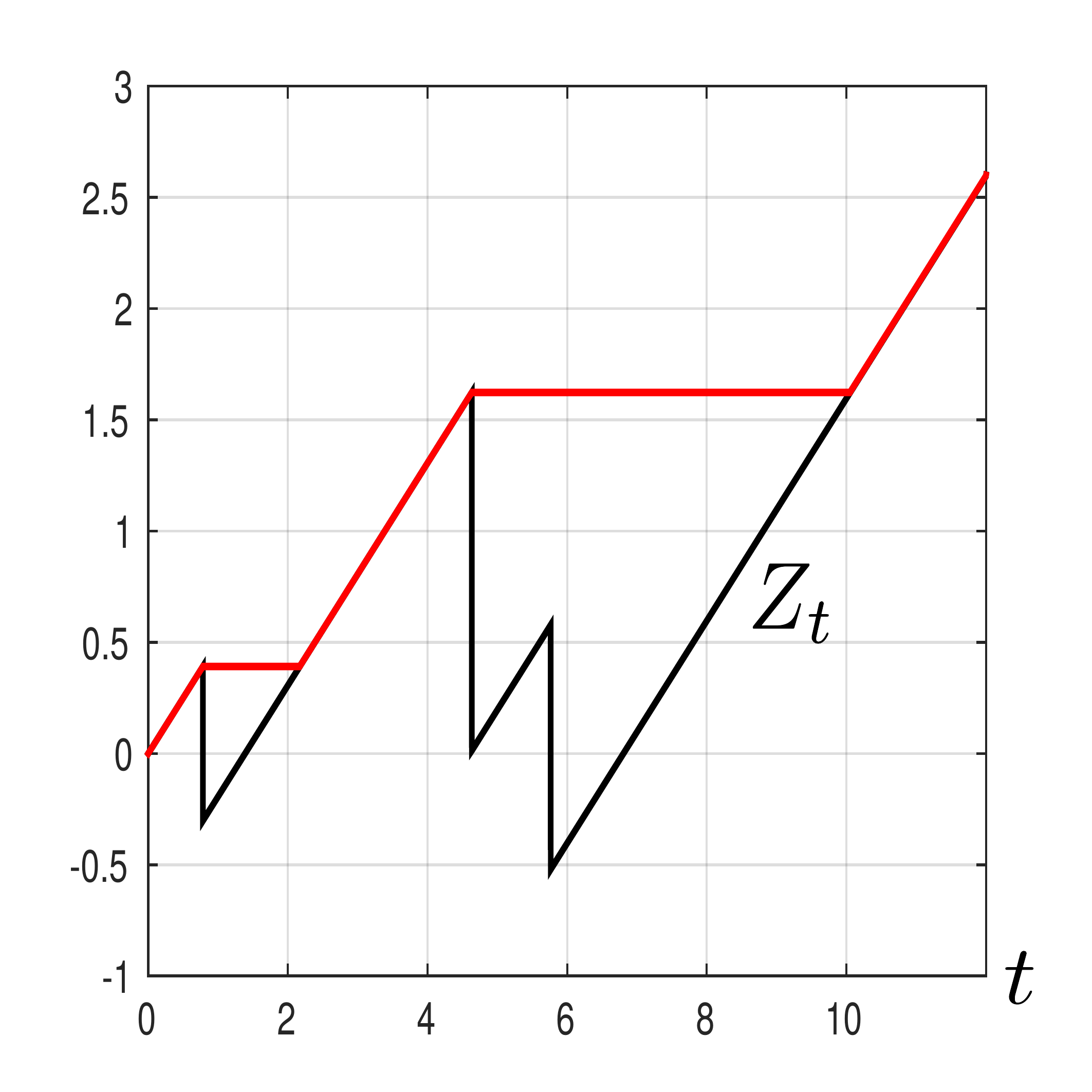}} 
\subfloat[][]{\label{f13}\includegraphics[height =5.5cm]{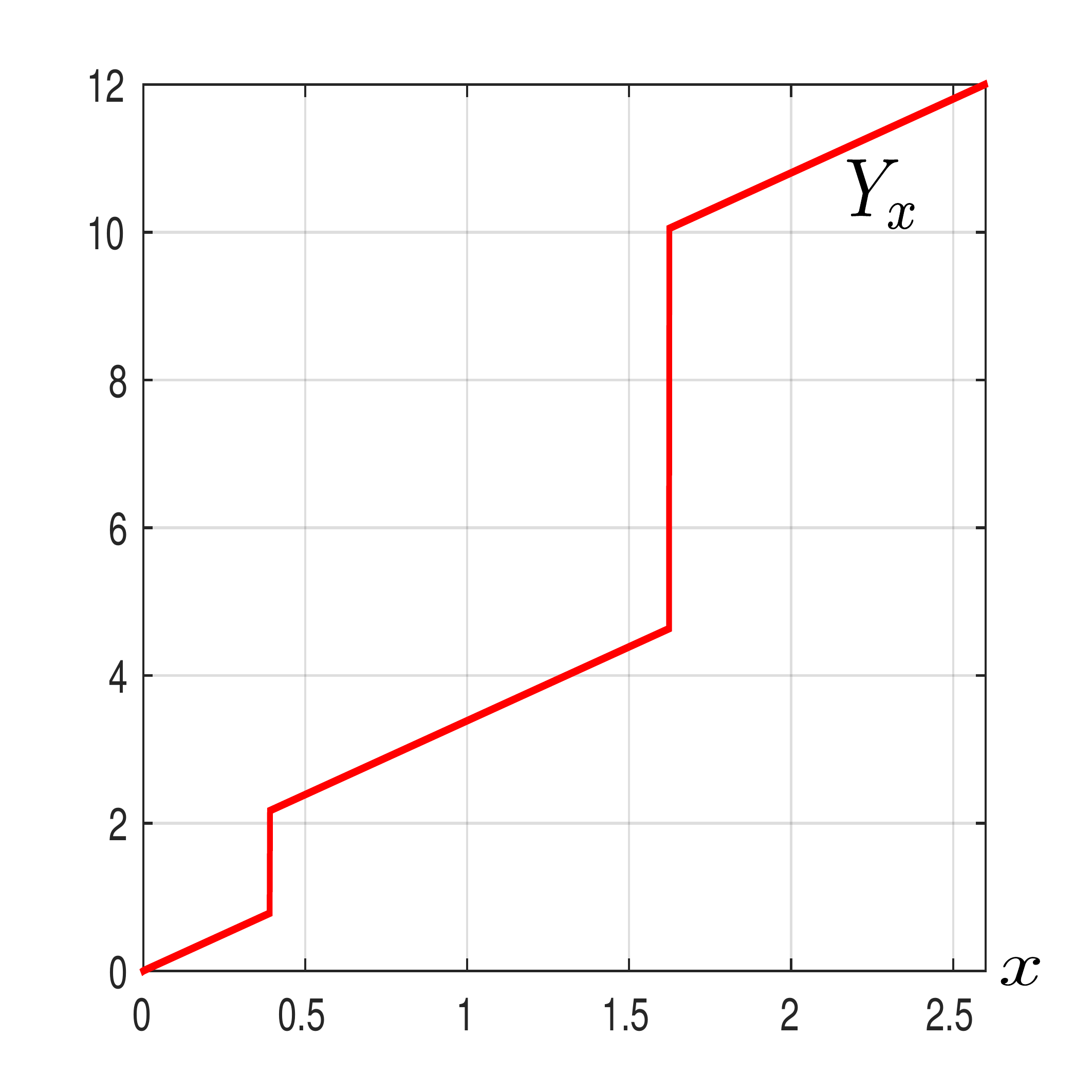}} 
\caption{
The relationship between subordinators $X$, $Y$ and the spectrally-negative process $Z$. Here $c=2$ and $0\le t \le 12$.
We have $\max(Z_t : \; 0\le t \le 12)=Z(12)\approx 2.6$, so the $x$-axis on graph $(c)$ has range $0\le x \le 2.6$. The red curve in figure (c) is obtained as a reflection of the red curve in figure (b) with respect to the diagonal line (the same transformation that we would use to find the graph of the inverse function).  
}
\label{fig1}
\end{figure}

Now, given a subordinator $X$, we fix $c>0$ and define a process $Z_t=t/c-X_t$. For $x\ge 0$ we introduce
\begin{equation}\label{def_Yx}
Y_x:=\inf\{t> 0: Z_t>x\}, 
\end{equation}
where we set $Y_x=+\infty$ on the event $\max\{ Z_t :  t\ge 0\} \le  x$.  The relationship between processes 
$X$, $Y$ and $Z$ can be seen on Figure \ref{fig1}. 
The process $Z_t=t/c-X_t$ is an example of a  {\it spectrally-negative} L\'evy process, and the random variables $Y_x$ are called the {\it first-passage times}, see Chapter 3 in \cite{Kyprianou}. 
It is known that the process $Y=\{Y_x\}_{x\ge 0}$ is a (possibly killed) subordinator, see   \cite{Kyprianou}[Corollary 3.14], 
thus there exists a Bernstein function $\phi_Y(z)$ such that for $w>0$ and $x>0$ we have 
\begin{equation}\label{Laplace_Yx}
\e[e^{-w Y_x}]=e^{-x \phi_Y(w)}.
\end{equation}
The Bernstein function $\Phi_Y(w)$ is known to satisfy the functional equation
\begin{equation}\label{functional_equation_Phi_Y}
 z/c - \phi_X(z)=w, \;\;\ w>0 \Longleftrightarrow z=\phi_Y(w), 
\end{equation}
see  \cite{Kyprianou}[Theorem 3.12]. Moreover, the distributions of $\{Y_x\}_{x\ge 0}$ and $\{X_t\}_{t\ge 0}$ are related through Kendall's identity  
(see  \cite{ECP1038}, \cite{Burridge}, \cite{Kendall} or  \cite{Kyprianou}[exercise 6.10]) 
\beq\label{Kendalls_identity}
\int_y^{\infty} \p(Y_x \le t) \frac{\d x}{x}=\int_0^t \p(Z_s > y) \frac{\d s}{s}, \;\;\; y>0, \;\;\; t>0.   
\eeq
The above identity will be the main ingredient in our proof of Theorem \ref{thm1}.

\subsection{Probabilistic proof of the case when $a(n)\ge 0$, $v>0$ and $w>0$}

We recall that the von Mangoldt function is defined via \eqref{def_Lambda} and we introduce a positive measure 
\begin{equation}\label{def_Pi}
\Pi(\d x)=\sum\limits_{n\ge 2} \frac{ \Lambda(n) a(n)}{\log(n)  n^{\sigma}} \delta_{\ln(n)}(\d x),
\end{equation}
where $\delta_y(\d x)$ denotes the Dirac measure concentrated at point $y$.
In the above formula $\sigma$ is chosen so that the Dirichlet series \eqref{Dirichlet_series_for_L} converges absolutely for $\re(s)\ge \sigma$. 
For $\re(s)\ge \sigma$ we have
\begin{align}\label{L_infinite_product}
L(s)&=\prod\limits_{p} (1-a(p)p^{-s})^{-1}=\exp(-\sum\limits_{p} \ln(1-a(p)p^{-s}))\\
\nonumber
&=
\exp\Big( \sum\limits_{p} \sum\limits_{k\ge 1} \frac{a(p)^k}{k p^{ks}}\Big)
=\exp\Big( \sum\limits_{n\ge 2} \frac{\Lambda(n)a(n) }{\log(n)n^s }\Big),
\end{align}
which implies that $\Pi$ is a finite measure of total mass $\Pi((0,\infty))=\ln(L(\sigma))$. 
Let us now consider a Bernstein function corresponding to the measure $\Pi$, that is
\begin{equation}\label{def_phi_X1}
\phi_X(z)=\int\limits_{(0,\infty)} (1-e^{-zx}) \Pi(\d x).
\end{equation}
From formulas \eqref{def_Pi} and \eqref{L_infinite_product} we see that 
\begin{equation}\label{def_phi_X}
\phi_X(z)=-\ln(L(\sigma+z))+\ln(L(\sigma)). 
\end{equation}
Let $X$ be a compound Poisson subordinator associated to the Bernstein function $\phi_X$. Comparing \eqref{Laplaceexponent} and \eqref{def_phi_X1} we see that the process $X$ has zero killing rate and zero linear drift, so it is a pure jump compound Poisson process. 
Due to \eqref{Laplace_X_t} and \eqref{def_phi_X} this process satisfies
\begin{equation}\label{X_t_Laplace}
\e[e^{-z X_t}]=e^{-t \phi_X(z)}=\Big( \frac{L(\sigma+z)}{L(\sigma)} \Big)^t. 
\end{equation}

We would like to point out that the above observations are not new: in the case $L(s)=\zeta(s)$ it was observed by Khintchine \cite{Khinchine} back in 1938 that $L(\sigma+\i z)/L(\sigma)$ is a characteristic function of an infinitely divisible distribution, and, more recently, the connections between more general L-functions and infinite divisibility were studied in \cite{donglin2001}.

Using binomial series we can easily find the measure $\p(X_t \in \d x)$. For $\re(s)\ge \sigma$ and $t>0$ we calculate 
\begin{equation}
\label{eqn_L_s_t}
L(s)^t=\prod\limits_{p} (1-a(p)p^{-s})^{-t}=
\prod\limits_{p} \sum\limits_{j \ge 1} \binom{j+t-1}{j} \frac{a(p)^j}{p^{js}}
=\sum\limits_{n\ge 1} d_t(n) \frac{a(n)}{n^s}.
\end{equation}
The above formula combined with \eqref{X_t_Laplace} show that for every $t>0$ the random variable $X_t$ is supported on the set
$\{\ln(n)\}_{n\in {\mathbb N}}$ and 
\begin{equation}\label{distribution_of_Xt}
\p(X_t=\ln(n))=L(\sigma)^{-t} d_t(n) \times \frac{a(n)}{n^{\sigma}}. 
\end{equation}

Now we fix $c>0$, denote $Z_t=t/c-X_t$ and we define the subordinator $Y$ as in \eqref{def_Yx}. 
The goal now is to use Kendall's identity to find the distribution of $Y_x$. The calculation that follows will be very similar to the one performed in the proof Proposition 3 in \cite{Burridge}. 

First of all, we claim that for every $x>0$ the random variable $Y_x$ has support on the set                   
$\{cx+c\ln(n)\}_{n\in {\mathbb N}}$.  This can be seen as follows. If the spectrally negative process $Z_t=t/c-X_t$ has no jumps before it hits level $x$, then $Y_x=cx$; if $Z$ has one jump of size $\ln(n_1)$ before it hits the level $x$, then $Y_x=cx+c\ln(n_1)$; if  $Z$ has two jumps of size
$\ln(n_1)$ and $\ln(n_2)$ before it hits level $x$, then $Y_x=cx+c\ln(n_1n_2)$, etc. Thus in order to describe the distribution of $Y_x$ it is enough to compute 
$$
p(n,x):=\p(Y_x=cx+c\ln(n)), \;\;\; n \in {\mathbb N}, \;  x>0.  
$$

The left-hand side of Kendall's identity \eqref{Kendalls_identity} can be written in the following way 
\begin{align}\label{lhs_Kendall}
\int_y^{\infty} \p(Y_x \le t) \frac{\d x}{x}=\int_y^{\infty} 
\sum\limits_{n\ge 1} {\mathbb I}_{\{cx+\ln(n)\le t\}} p(n,x) \frac{\d x}{x}
=\sum\limits_{1\le n \le \exp(t/c-y)}
\int_y^{t/c-\ln(n)} p(n,x) \frac{\d x}{x}. 
\end{align}
And the right-hand side of \eqref{Kendalls_identity} is transformed into
\begin{align}\label{rhs_Kendall}
\nonumber
\int_0^t \p(X_s < s/c- y) \frac{\d s}{s}&=\int_0^t \sum\limits_{1\le n < \exp(s/c - y)}
\p(X_s=\ln(n)) \frac{\d s}{s} \\
&= 
\sum\limits_{1\le n < \exp(t/c-y)} \int_{cy+c\ln(n)}^t \p(X_s=\ln(n)) \frac{\d s}{s}\\
\nonumber
&= 
\sum\limits_{1\le n < \exp(ct-y)} \int_{y}^{t/c-\ln(n)} \p(X_{cu+c\ln(n)}=\ln(n)) \frac{\d u}{u+\ln(n)},
\end{align}
where in the last step we have changed the variable of integration $s=cu+c\ln(n)$. 
Using Kendall's identity \eqref{Kendalls_identity} and comparing the expressions in the right-hand side in 
\eqref{lhs_Kendall} and \eqref{rhs_Kendall} we conclude that 
\begin{equation*}
p(n,x)=\p(X_{cx+c\ln(n)}=\ln(n)) \frac{x}{x+\ln(n)}.
\end{equation*}
Applying \eqref{distribution_of_Xt} to the above formula we see that for $n\in {\mathbb N}$
\begin{align*}
\p(Y_x=cx+c\ln(n))&=p(n,x)=L(\sigma)^{-cx-c\ln(n)} d_{cx+c\ln(n)}(n) \times \frac{a(n)}{n^{\sigma}} \times \frac{x}{x+\ln(n)}\\
&=L(\sigma)^{-cx} d_{cx+c\ln(n)}(n) \times \frac{a(n)}{n^{\sigma+c\ln(L(\sigma))}} \times \frac{x}{x+\ln(n)}.
\end{align*}
Next, combining the above result with \eqref{Laplace_Yx} we conclude that for $x>0$ and $w>0$
\begin{align*}
\e[e^{-w Y_x}]&=\sum\limits_{n\ge 1} e^{-w (cx+c\ln(n))} \p( Y_x=cx+c\ln(n))\\
&=   \sum\limits_{n\ge 1}
e^{-cwx} n^{-cw}
 L(\sigma)^{-cx} d_{cx+c\ln(n)}(n) \times \frac{a(n)}{n^{\sigma+c\ln(L(\sigma))}} \times \frac{x}{x+\ln(n)}=
 e^{-x \phi_Y(w)}. 
\end{align*}
We use our definition of $\tilde d_t(n)=t^{-1} d_t(n)$, rearrange the terms in the above identity and rewrite it in the form
\begin{equation}\label{eqn_identity1}
 \sum\limits_{n\ge 2}
 \tilde d_{cx+c\ln(n)}(n) \times \frac{a(n)}{n^{\sigma+c(w+\ln(L(\sigma)))}}=
 \frac{1}{cx} \Big( e^{x(cw-\phi_Y(w)+c\ln(L(\sigma)))}-1 \Big).
\end{equation}
We emphasize that formula \eqref{eqn_identity1} is valid for all $x > 0$ and $w>0$ and that $z=\phi_Y(w)$ is the solution to the equation   
\begin{equation}\label{functional_eqn_L}
z/c+\ln(L(\sigma+z))-\ln(L(\sigma))=w, 
\end{equation}
see \eqref{functional_equation_Phi_Y} and \eqref{def_phi_X}. 
 
Let us introduce a new function 
\begin{equation}\label{def_f_s_c}
f(s,c)=c^{-1}\times (s-\phi_Y((s-\sigma)/c-\ln(L(\sigma)))-\sigma), \;\;\; s>\sigma+c\ln(L(\sigma)).  
\end{equation}
Using \eqref{functional_eqn_L} we check that $f$ satisfies the equation $L(s-cf(s,c))=\exp(f(s,c))$ and 
formula \eqref{eqn_identity1} can be rewritten as 
 \begin{equation}\label{eqn_identity2}
 \sum\limits_{n\ge 2}
 \tilde d_{cx+c\ln(n)}(n) \times \frac{a(n)}{n^{s}}=
 \frac{1}{cx} \Big( e^{c x f(s,c)}-1 \Big), \;\;\; s>\sigma+c\ln(L(\sigma)).
\end{equation}
This ends the proof of  \eqref{inverse_eqn} and  \eqref{eqn_exp_v_f_s_w_}. Formula \eqref{def_f_s_w} is derived 
 by taking the limit in \eqref{eqn_identity2} as $x \to 0^+$. 
For $s>\sigma+c\ln(L(\sigma))$, the lower bound $0<f(s,w)$ follows from \eqref{def_f_s_w}, and the upper bound $f(s,w)<\gamma$ can be easily established from \eqref{def_f_s_c} and the fact that   
$$
\phi_Y(w)=cw+c\phi_X(\phi_Y(w))>cw,
$$
which is a consequence of \eqref{functional_equation_Phi_Y}. Thus $0<f(s,w)<\gamma$ for $s>\sigma+c\ln(L(\sigma))$, and since 
$|f(s,w)|\le f(\re(s),w)$ the result holds for complex $s$ in the half-plane $\re(s)>\sigma+c\ln(L(\sigma))$ as well. 

Thus we have proved all statements of Theorem \ref{thm1} in the case $v=cx>0$, $w=c>0$, and $a(n)\ge 0$.

\subsection{Proving the general result via analytic continuation}

So far we have proved that Theorem \ref{thm1} holds for $v>0$, $w>0$ and $a(n)\ge 0$. Our first goal is to extend this result to complex values of $v$ and $w$. The key observation is that $\tilde d_t(n)$ is a polynomial of $t$ whose roots are non-positive integers. 
Writing this polynomial as a product of linear factors and applying the inequality $|q+t|\le q+|t|$ (with $q>0$ and $t\in \C$) to each linear factor, we deduce the upper bound  
\begin{equation}\label{estimate_tau}
|\tilde d_t(n)|\le \tilde d_{|t|}(n), \;\;\; n\ge 2, \;\;\; t\in \C. 
\end{equation}  
Therefore, if the series \eqref{def_f_s_w} converges for some $w=\rho>0$ and $s=\sigma+\gamma \rho$, it will converge uniformly for $(s,w) \in D_{\sigma+\gamma \rho,\rho}$. From this fact we see that the function $f(s,w)$ is an analytic function of two variables $(s,w) \in D_{\sigma+\gamma \rho,\rho}$. Moreover, the inequality \eqref{estimate_tau} implies that $|f(s,w)|\le f(\re(s),|w|)<\gamma$ for 
$(s,w) \in D_{\sigma+\gamma \rho,\rho}$.
Since $L(s)$ is analytic in $\re(s)>\sigma$ and the function $f(s,w)$ satisfies the functional equation \eqref{inverse_eqn} for $w\in (0,\rho)$, we conclude by analytic continuation that the same equation must hold for all $(s,w) \in D_{\sigma+\gamma \rho,\rho}$. 
Thus we have extended Theorem \ref{thm1} to allow for complex values of $w$. 
To prove \eqref{eqn_exp_v_f_s_w_} in the general case when $v$ is complex, we use the same approach and an analytic continuation in $v$. 

Our goal now is to remove the remaining restriction -- the condition that $a(n)\ge 0$.
Let us denote $i$-th prime number by $p_i$ (so that $p_1=2$, $p_2=3$, $p_3=5$, etc.). 
Consider ${\mathbf u}=(u_1,u_2,\dots,u_k)\in \C^k$ and a Dirichlet L-function
\begin{equation}
\label{eqn_finite_Euler}
L(s;{\mathbf u})=\prod\limits_{i=1}^k (1-u_i p_i^{-s})^{-1}.
\end{equation}
Denote by $a(n;{\mathbf u})$ the corresponding completely multiplicative function, that is 
\begin{equation}\label{eqn_a_n_u}
a(n;{\mathbf u})=
u_1^{l_1} \dots u_k^{l_k} \;\; {\textnormal{ if }} \; n=p_1^{l_1} \dots p_k^{l_k}
\end{equation}
and $a(n;{\mathbf u})=0$ otherwise. 
Let us now fix 
${\mathbf u}\in \C^k$ and denote
$$
B({\mathbf u}):=\{ {\mathbf x} \in \C^k \; : \; |x_i|\le |u_i|\} 
{\textnormal{  and  }} C({\mathbf u}):=\{ {\mathbf x} \in \r^k \; : \; 0\le x_i \le |u_i|\}. 
$$

For ${\mathbf x}\in C({\mathbf u})$ the completely multiplicative function $a(n;{\mathbf x})$ is non-negative, thus Theorem \ref{thm1} holds for $L(s;{\mathbf x})$. 
Consider a function $L(s;|{\mathbf u}|)$, where $|{\mathbf u}|:=(|u_1|,\dots,|u_n|)$. There exists $\sigma$ such that the Dirichlet series for this function converges absolutely for $\re(s)\ge \sigma$. Note that for any  
 ${\mathbf x}\in C({\mathbf v})$ the Dirichlet series for $L(s;{\mathbf x})$ also converges absolutely 
in $\re(s)\ge \sigma$, since $|a(n;{\mathbf x})|\le a(n;|{\mathbf u}|)$. 

Let us denote 
\begin{equation}\label{def_gamma1}
\gamma=\gamma(|{\mathbf u}|)=\ln \Big(\sum\limits_{n\ge 1} \frac{a(n;|{\mathbf u}|)}{n^{\sigma}}\Big). 
\end{equation}
Applying Theorem \ref{thm1} to each Dirichlet L-function $L(s;{\mathbf x})$ we conclude that for any $\rho>0$ and any ${\mathbf x} \in C({\mathbf u})$ the following results hold:
\begin{itemize}
\item[(i)] The series
\begin{equation}\label{def_f_s_w1}
f(s,w;{\mathbf x}):=\sum\limits_{n\ge 2} \tilde d_{w\ln(n)}(n)\times \frac{a(n;{\mathbf x})}{n^s} 
\end{equation}
converges absolutely and uniformly for all $(s,w) \in D_{\sigma+\gamma\rho,\rho}$ and satisfies $|f(s,w;{\mathbf x})|<\gamma$ in this region; 
\item[(ii)] 
The function $f(s,w;{\mathbf x})$ solves the functional equation
\begin{equation}\label{inverse_eqn1}
L(s-wf(s,w;{\mathbf x});{\mathbf x})=\exp(f(s,w;{\mathbf x})), \;\;\; (s,w)  \in D_{\sigma+\gamma\rho,\rho}; 
\end{equation}
\item[(iii)]
For any $v \in \C \setminus\{0\}$ and $(s,w) \in D_{\sigma+\gamma\rho,\rho}$ the following identity is true
\begin{equation}\label{eqn_exp_v_f_s_w_1}
1+v\sum\limits_{n\ge 2}\tilde d_{v+w\ln(n)}(n) \times \frac{a(n;{\mathbf x})}{n^s}= \exp(vf(s,w;{\mathbf x}).
\end{equation}
\end{itemize}

Let us now fix values of $v \in \C$ and $(s,w) \in D_{\sigma+\gamma \rho, \rho}$. 
Note that the function ${\mathbf x} \mapsto L(s;{\mathbf x})$ is analytic in  the interior of the set $B({\mathbf u})$ and continuous in $B({\mathbf u})$.
Next, formula \eqref{eqn_a_n_u} shows that $a(n;{\mathbf x})$ are monomials in variables $x_1, x_2, \dots, x_k$,
thus formula \eqref{def_f_s_w1}, which defines the function  $f(s,w;{\mathbf x})$, can be viewed as a Taylor series for the function ${\mathbf x} \mapsto f(s,w;{\mathbf x})$. Since $|f(s,w;{\mathbf x})|\le f(\re(s),|w|;|{\mathbf x}|)< \gamma$, this Taylor series converges uniformly for all ${\mathbf x} \in B({\mathbf u})$, so that ${\mathbf x} \mapsto f(s,w;{\mathbf x})$ is an analytic function in the interior of the set $B({\mathbf u})$, and it is continuous in $B({\mathbf u})$. 
Using these results and analytic continuation, we can extend the functional equation \eqref{inverse_eqn1} for all ${\mathbf x} \in B({\mathbf u})$ (since we have already established that \eqref{inverse_eqn1} holds true for 
${\mathbf x} \in C({\mathbf u})$). 

The same argument applies to the identity \eqref{eqn_exp_v_f_s_w_1}. The right-hand side is an analytic function of ${\mathbf x}$ in the interior of the set $B({\mathbf u})$, and it is continuous in $B({\mathbf u})$. The left-hand side is a power series in ${\mathbf x}$, convergent uniformly in $B({\mathbf u})$. By analytic continuation, since the identity \eqref{eqn_exp_v_f_s_w_1} holds for ${\mathbf x} \in C({\mathbf u})$, it must hold everywhere in 
$B({\mathbf u})$. 

Thus we have shown that formulas \eqref{inverse_eqn1} and \eqref{eqn_exp_v_f_s_w_1} hold for all
$v\in\C$, $(s,w) \in D_{\sigma+\gamma \rho, \rho}$ and ${\mathbf x} \in 
B({\mathbf u})$. In particular, they must hold for ${\mathbf x}={\mathbf u}$. This ends the proof of Theorem \ref{thm1} for L-functions of the form \eqref{eqn_finite_Euler}. 

Finally, consider a general Dirichlet L-function $L(s)$ defined via \eqref{Dirichlet_series_for_L}. 
Assume that the Dirichlet series for $L(s)$ converges absolutely for $\re(s)\ge \sigma$. Let us define 
\begin{equation}\label{def_L_k}
L_k(s)=\prod\limits_{i=1}^k (1-a(p_i)p_i^{-s})^{-1}, \;\;\; k\ge 1,
\end{equation}
and 
\begin{equation}\label{def_tilde_L}
\tilde L(s)=\prod\limits_{i=1}^{\infty} (1-|a(p_i)|p_i^{-s})^{-1}.
\end{equation}
It is clear that the Dirichlet series for all L-functions $L_k(s)$ and $\tilde L(s)$ converge absolutely when $\re(s)\ge \sigma$. 
Let us denote by $a_k(n)$ the completely multiplicative function corresponding to the L-function $L_k(s)$ and 
let $\gamma$ be defined via \eqref{def_gamma}. 
We have proved already that Theorem \ref{thm1} holds true for $L_k(s)$ and $\tilde L(s)$. Thus the following results hold true. For any $\rho>0$:
\begin{itemize}
\item[(i)] The series
\begin{equation}\label{def_f_s_w2}
\tilde f(s,w):=\sum\limits_{n\ge 2} \tilde d_{w\ln(n)}(n)\times \frac{|a(n)|}{n^s} 
\end{equation}
converges absolutely and uniformly for all $(s,w) \in D_{\sigma+\gamma\rho,\rho}$ and satisfies $|\tilde f(s,w)|<\gamma$ in this region; 
\item[(ii)] For each $k\ge 1$, the series
\begin{equation}\label{def_f_s_w3}
f_k(s,w):=\sum\limits_{n\ge 2} \tilde d_{w\ln(n)}(n)\times \frac{a_k(n)}{n^s} 
\end{equation}
converges absolutely and uniformly for all $(s,w) \in D_{\sigma+\gamma\rho,\rho}$ and satisfies $|f_k(s,w)|<\gamma$ in this region; 
\item[(iii)] 
For each $k\ge 1$, the functions $f_k(s,w)$ solve the functional equation
\begin{equation}\label{inverse_eqn2}
L_k(s-wf_k(s,w))=\exp(f_k(s,w)), \;\;\; (s,w)  \in D_{\sigma+\gamma\rho,\rho}; 
\end{equation}
\item[(iv)]
For any $v \in \C \setminus\{0\}$ and $(s,w) \in D_{\sigma+\gamma\rho,\rho}$ the following identity is true
\begin{equation}\label{eqn_exp_v_f_s_w_2}
1+v\sum\limits_{n\ge 2}\tilde d_{v+w\ln(n)}(n) \times \frac{|a(n)|}{n^s}= \exp(v\tilde f(s,w));
\end{equation}
\item[(v)]
For any $k\ge 1$, $v \in \C \setminus\{0\}$ and $(s,w) \in D_{\sigma+\gamma\rho,\rho}$ the following identity is true
\begin{equation}\label{eqn_exp_v_f_s_w_3}
1+v\sum\limits_{n\ge 2}\tilde d_{v+w\ln(n)}(n) \times \frac{a_k(n)}{n^s}= \exp(v f_k(s,w)).
\end{equation}
\end{itemize}

Note that the function $f(s,w)$ in \eqref{def_f_s_w} is well-defined, since the series converges absolutely for all
$(s,w) \in D_{\sigma+\gamma\rho,\rho}$, due to the absolute convergence of \eqref{def_f_s_w2}. 
It is clear from our definition that for each $n\in {\mathbb N}$ we have $a_k(n) \to a(n)$ as $k\to +\infty$
and that $|a_k(n)|\le |a(n)|$ for all $k, n \in {\mathbb N}$. 
Thus, we can use the Dominated Convergence Theorem, the fact that 
the series in \eqref{def_f_s_w2} converges absolutely and, by taking the limit as $k\to +\infty$ in \eqref{def_f_s_w3},
 we conclude that for all $(s,w) \in D_{\sigma+\gamma\rho,\rho}$ it is true that $f_k(s,w) \to f(s,w)$ as $k\to +\infty$. Since the functions $L_k(s)$ converge to $L(s)$ uniformly in the half-plane $\re(s) \ge \sigma$, we can take the limit as $k\to+\infty$ in \eqref{inverse_eqn2} and conclude that the functional equation 
\eqref{inverse_eqn} holds for all $(s,w) \in D_{\sigma+\gamma\rho,\rho}$. 
Finally, formula \eqref{eqn_exp_v_f_s_w_} can be established by taking the limit as $k\to +\infty$ in \eqref{eqn_exp_v_f_s_w_3} and applying the Dominated Convergence Theorem (with the help of the absolute convergence in \eqref{eqn_exp_v_f_s_w_2}). 

This ends the proof of Theorem \ref{thm1}.

\section{Proof of Corollary \ref{cor1}}\label{section_cor_proofs}

Assume that  $(s,w) \in D_{\sigma+\gamma \rho, \rho}$ for some $\rho >0$, which is equivalent to saying that 
$\re(s)\ge \sigma + \gamma |w|$. Denote $\alpha:=s-w f(s,w)$. Since $|f(s,w)|<\gamma$, $\re(s)\ge \sigma+\gamma \rho$ and $|w|\le \rho$,  we have 
$\re(\alpha)>\sigma$. 
Identity \eqref{inverse_eqn} tells us that $L(\alpha)=\exp(f(s,w))$, so that $f(s,w)=\ln (L(\alpha))=:\beta$. Moreover, 
from equation $\alpha=s-w f(s,w)$ we express  $s=\alpha+w f(s,w)=\alpha + \beta w$ and then the equation $f(s,w)=\beta$ gives us 
\begin{equation}\label{cor1_proof1}
f(\alpha+\beta w,w)=\beta. 
\end{equation}
We emphasize that \eqref{cor1_proof1} holds for all $(\alpha,w) \in \C^2$ such that $\alpha=s-wf(s,w)$ for some $(s,w) \in \C^2$ satisfying 
$\re(s)\ge \sigma+\gamma |w|$. By analytic continuation we can extend \eqref{cor1_proof1} to all $(\alpha,w)$ such that 
\begin{equation}\label{cor1_proof2}
\re(\alpha)\ge \sigma \;\;\; {\textnormal{ and }} \;\;\; \re(\alpha + \beta w)\ge \sigma + \gamma |w|, 
\end{equation}
where $\beta=\ln(L(\alpha))$. 
 
The last step is to prove that condition $\re(\alpha)\ge \sigma + 2\gamma |w|$ implies \eqref{cor1_proof2}. This follows from the following sequence of inequalities
$$
\re(\alpha + \beta w)\ge \re(\alpha) - |\re(\beta w)|\ge \sigma + 2 \gamma |w| - |\beta| |w|\ge \sigma + \gamma |w|,
$$
where in the last step we have used the fact that $|\beta|=|\ln(L(\alpha))|\le \gamma$, which follows from  
\eqref{def_gamma} and \eqref{eqn_f_s_0}.  
\qed


\section*{Acknowledgements}
The author would like to thank Aleksandar Ivi\'c for comments and for pointing out relevant literature. The research is supported by the Natural Sciences and Engineering Research Council of Canada.





\end{document}